# Estimation of mean using Dual-to-Ratio and Difference-type estimators Under Measurement Error model

Viplav Kumar Singh and Rajesh Singh*

Department of Statistics, Banaras Hindu University,

Varanasi -221005, India

*Corresponding Author

**Abstract**

In sample survey, when data is collected, it is assumed that whatever is reported by respondent is correct. However, given the issues of prestige bias, personal respect, respondent's self reported data often produces over-or-under estimated values from true value. This causes measurement error to be present in sample values. In support of this study, we have considered some precise classes using dual under measurement error model. The expressions for the bias (B) and the mean square errors (MSE) of proposed classes have been derived and compared with, the mean per unit estimator, the Srivenkataramana (1980) estimator and Sharma and Tailor (2010) estimator.

## 1. Introduction

In past few decade's, Statisticians have paid their attention towards the problem of estimation of slope parameters in the presence of measurement errors. Basically, measurement error may be characterized as the difference between the value of a variable provided by the respondent and the true value of the same variable. The total survey error of a statistics with measurement error has both fixed bias error and variable error (variance) over repeated trails of the survey [see Sukhatme et al.(1984); Cochran (2005)].Figure 1 illustrates the concept of measurement error:

**Figure 1**

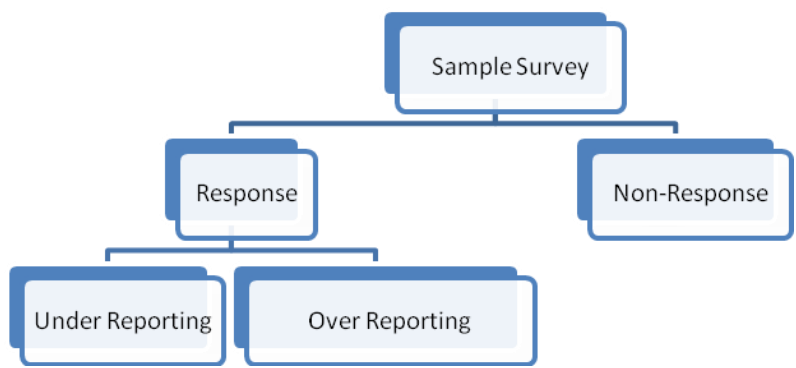

**Remark** In figure 1, under reporting and over reporting cause meaearment error.

Incompleteness in survey data may arrised due to: incorrect response or non-response. Measurment bias provides a systematic pattern in the difference between the respondents answer to a question and the correct answer. For example, 1. The survey interwiever asking about deaths were poorly trained asking about deaths poorly trained and included deaths which occurred before the time period of interest. This would lead to an overestimate of the mortality rate because deaths which should not be included are included. 2. One survey team's portable machine to measure haemoglobin malfunctioned and was not checked, as should be done every day. It measured everyone's haemoglobin as 0.3 g/L too high. This would lead to an underestimate of the prevalence of anaemia because the readings would overestimate the haemoglobin for everyone measured by that team.

Further, measurment variance reflects random variation in answers provided to an interviewer while asking the same question, that is, often the same respondent provides different answers to the same question when asked repeatedly. Several methods are available in the survey sampling literature to handle non-response, including the revisit method, imputation methods, auxiliary sources utilization method and the neighbouring units manipulation methods, however, when a respondent provides incorrect information regarding a variable, additional techniques are required. This study considers this aspect and deals with mean estimation under measurement error.

Many researchers have paid their attention towards the problem of estimation of population parameters in the presence of measurement errors. Starting form Cochran (1968), who had studied the effect of measurement error on the data analysis. Shalabh (1997), Manisha and Singh (2001), Kumar et al.(2011) and Shukla et al.(2012) have addressed the problem of estimation of mean using information on auxiliary variable in the presence of measurement errors. Later, Srivastava and Shalabh (2001), Manisha and Singh (2002), Allen et al. (2003) and Singh and Karpe (2008, 2009) and others have made some more contribution on measurement errors. However, no effort has been made to estimate the finite population mean using dual-estimator in the presence of measurement error. This motivation led us to consider the problem of estimation of finite population mean using dual-to-ratio and difference type estimators in the presence of measurement error. In this paper, we adapted Srivastava (1971), Singh and Solanki (2012) and Sabbir and Yaab (2003) estimator and use it for estimating mean in the presence of measurement error. Expressions for the biases and mean square errors of adapted estimators have been derived up to the first order of approximation. An empirical study is also carried out to demonstrate the superiority of the adapted estimators over existing one.

## 2. Notations and Expectations

Let us consider a finite population $U = [U_1, U_2, ......., U_N]$ of size N. Let Y and X be the study and auxiliary variate, respectively. Suppose that a sample of size n is drawn using simple random sampling without replacement. It is assumed that $y_i$ and $x_i$ for the ith sampling unit are recorded with measurement error instead of their true values $X_i$ and $Y_i$ as

$$b_{Yi} = y_i - Y_i, \qquad (1)$$

and

$$b_{Xi} = x_i - X_i \qquad (2)$$

where $(b_{Yi}, b_{Xi})$ are the associated measurement errors which are assumed to be stochastic with mean zero and variances $S_{b_Y}^2$ and $S_{b_X}^2$ respectively. For simplicity, we assumed that $b_{Yi}$ and $b_{Xi}$ are uncorrelated although $X_i$'s and $Y_i$'s are correlated. We further assumed that the measurement errors are independent of true values of the variables.

Let $(\mu_x, \mu_y)$ and $(S_x^2, S_y^2)$ be the population means and variances of the characteristics (X, Y), respectively. Further, let $\rho$ be the population correlation coefficient between Y and X. Let $\bar{x} = \frac{1}{n}\sum_{i=1}^{n} x_i$, $\bar{y} = \frac{1}{n}\sum_{i=1}^{n} y_i$ be the unbiased estimators of the population means $\mu_x$ and $\mu_y$, respectively. Also, Let $C_Y$ and $C_X$ be the population co-efficient of variation for the variable Y and X respectively. We further assumed that the mean of the study variable Y is unknown and mean of auxiliary variable X is known.

In order to derive the bias and mean square error of the adapted estimators in the presence of measurement error, let us define the following notations.

Let

$$w_Y = \frac{1}{\sqrt{n}}\sum_{i=1}^{n}(Y_i - \bar{Y}), \tag{3}$$

$$w_X = \frac{1}{\sqrt{n}}\sum_{i=1}^{n}(X_i - \bar{X}), \tag{4}$$

$$w_{d_Y} = \frac{1}{\sqrt{n}}\sum_{i=1}^{n} d_{Y_i}, \tag{5}$$

and

$$w_{d_X} = \frac{1}{\sqrt{n}}\sum_{i=1}^{n} d_{X_i} \tag{6}$$

Adding (3) and (5), we have

$$w_Y + w_{d_Y} = \frac{1}{\sqrt{n}}\sum_{i=1}^{n}(Y_i - \bar{Y}) + \frac{1}{\sqrt{n}}\sum_{i=1}^{n} d_{Y_i} \tag{7}$$

Using (1) and (7) and simplifying, we get

$$\bar{y} = \bar{Y} + \frac{1}{\sqrt{n}}(w_Y + w_{d_Y}) = \bar{Y} + \kappa_Y \tag{8}$$

Similarly from (2), (4) and (6), we have

$$\bar{x} = \bar{X} + \frac{1}{\sqrt{n}}(w_X + w_{d_X}) = \bar{X} + \kappa_X \tag{9}$$

Furthermore

$$\left.\begin{array}{l} E(\kappa_Y) = \gamma(S_Y^2 + S_{d_Y}^2) = r_0 \\ E(\kappa_X) = \gamma(S_X^2 + S_{d_X}^2) = r_1 \\ E(\kappa_Y \kappa_X) = \gamma \rho S_Y S_X = r_{01} \end{array}\right\} \text{(say)} \tag{10}$$

where $\gamma = \left(\frac{1}{n} - \frac{1}{N}\right)$

## 3. Existing and Adapted estimators
In section 3.1, we have given some well known existing estimators in literature with their properties. Similarly, section 3.2 contains the adapted estimators with their properties.

### 3.1 Existing Estimators and their Properties

**3.1.1** The Mean per unit estimator is given by

$$\bar{y} = \frac{1}{n}\sum_{i=1}^{n} y_i \tag{11}$$

Estimator $\bar{y}$ is unbiased with variance, given as

$$V(\bar{y}) = \gamma \bar{Y}^2 \left( C_Y^2 + \frac{S_{d_Y}^2}{\bar{Y}^2} \right) \tag{12}$$

where $\gamma$ is already defined in section 2.

### 3.1.2 Srivenkataramana (1980) Estimator

Srivenkataramana (1980) suggested a dual-to-ratio type estimator as

$$e_1 = \bar{y}\frac{\bar{x}^*}{\bar{X}}, \tag{13}$$

where

$$\bar{x}^* = \frac{(N\bar{X} - n\bar{x})}{(N-n)}. \tag{14}$$

The MSE of $e_1$ up to the first order of approximation, is given as

$$MSE(e_1) = \gamma \bar{Y}^2 \left[ C_Y^2 + n_1^2 C_X^2 - 2n_1 \rho C_Y C_X \right] \tag{15}$$

where

$$n_1 = \frac{n}{(N-n)}.$$

If measurement errors are also taken into account, estimator $e_1$ becomes

$$e_1^m = \bar{y}\frac{\bar{x}^{**}}{\mu_x} \tag{16}$$

And their respective mean square error is given by

$$MSE(e_1^m) = \gamma \bar{Y}^2 \left[ C_Y^2 + n_1^2 C_X^2 - 2n_1 \rho C_Y C_X \right] + \gamma \left[ S_{d_Y}^2 + n_1^2 R^2 S_{d_X}^2 \right] \tag{17}$$

Here the second part of equation (17) is the contribution of measurement error to the mean square error of $e_1^m$.

where

$$\bar{x}^{**} = \frac{(N\mu_x - n\bar{x})}{(N-n)} \text{ and } R = \frac{\bar{Y}}{\mu_x}. \tag{18}$$

### 3.1.3 Sharma and Tailor (2010) Estimator

Sharma and Tailor (2010) suggested the following ratio-cum-dual to ratio estimator by taking the linear combination of classical ratio estimator and dual to ratio estimator.

$$e_2 = \bar{y}\left[\alpha \frac{\bar{X}}{\bar{x}} + (1-\alpha)\frac{\bar{x}^*}{\bar{X}}\right] \tag{19}$$

The MSE of $e_2$ can be expressed as

$$\text{MSE}(e_2) = \bar{Y}^2\left[1 + \alpha^2 A_1 + (1-\alpha)^2 A_2 - 2\alpha A_3 - 2(1-\alpha)A_4 + 2\alpha(1-\alpha)A_5\right] \tag{20}$$

Also, the minimum MSE of $e_2$ is obtained for optimum value of $\alpha$, given as

$$\alpha(\text{opt}) = \frac{(A_2 + A_3 - A_4 - A_4)}{(A_1 + A_2 - 2A_5)} = \alpha^\oplus \tag{21}$$

where

$$A_1 = 1 + \gamma\{C_y^2 + 3C_x^2 - 4\rho C_y C_x\}$$
$$A_2 = 1 + \gamma\{C_y^2 + n_1^2 C_x^2 - 4n_1\rho C_y C_x\}$$
$$A_3 = 1 + \gamma\{C_x^2 - \rho C_y C_x\}$$
$$A_4 = \{1 - n_1\gamma\rho C_y C_x\}$$
$$A_5 = 1 + \gamma\{C_y^2 + C_x^2(1+n_1) - 2\rho C_y C_x(1+n_1)\}$$

If measurement errors are taken into account, estimator $e_2$ becomes

$$e_2^m = \bar{y}\left[\alpha'\frac{\mu_x}{\bar{x}} + (1-\alpha')\frac{\bar{x}^{**}}{\mu_x}\right] \tag{22}$$

The MSE of estimator $e_2^m$ is given by

$$\text{MSE}(e_2^m) = \left[\bar{Y}^2 + \alpha'^2 B_1 + (1-\alpha')^2 B_2 - 2\alpha' B_3 - 2(1-\alpha')B_4 + 2\alpha(1-\alpha')B_5\right] \tag{23}$$

Also, the minimum MSE of $e_2^m$ is obtained for optimum value of $\alpha'$, given as

$$\alpha'(\text{opt}) = \frac{(B_2 + B_3 - B_4 - B_4)}{(B_1 + B_2 - 2B_5)} = \alpha'^\oplus \tag{24}$$

where

$$B_1 = \bar{Y}^2 + r_0 + 3R^2 r_1 - 4Rr_{01}$$
$$B_2 = \bar{Y}^2 + r_0 + n_1^2 R^2 r_1 - 4n_1 Rr_{01}$$
$$B_3 = \bar{Y}^2 + R^2 r_1 - Rr_{01}$$
$$B_4 = \bar{Y}^2 - n_1 Rr_{01}$$
$$B_5 = \bar{Y}^2 + r_0 + (1+n_1)R^2 r_1 - 2n_1 Rr_{01}(n_1+1)$$

Putting these values in (23) and (24), we have the min. MSE of $e_2^m$ for optimum value of $\alpha'$, respectively.

## 3.2 Adapted Estimators and their Properties

We have adapted Srivastava (1971), Sabbir and Yaab (2003) and Singh and Solanki (2012) estimators for estimating mean in the presence of measurement errors as follows

### 3.2.1 Wider class of estimators

Motivated by Srivastava (1971), we consider the following class of estimators using dual transformation in the presence of measurement error given as

$$\hat{\bar{Y}}_1 = g(\bar{y}, u^{**}) \tag{25}$$

where, $u^{**} = \left(\dfrac{\bar{x}^{**}}{\mu_x}\right)$ and $g(\bar{y}, u^{**})$ is a function of $\bar{y}$ and $u^{**}$ and satisfies the following regularity conditions

(i) The point $(\bar{y}, u^{**})$ assumes the value in the closed convex subset $R_2$ of two dimensional real space containing the point $(\bar{Y}, 1)$.

(ii) The function $g(\bar{y}, u^{**})$ is continuos and bounded in $R_2$.

(iii) $g(\bar{Y}, 1) = \bar{Y}$ and $G_0 = \dfrac{\partial g}{\partial \bar{y}} = 1$. Also, the first, second order derivatives of $g(\bar{y}, u^{**})$ exists and are continuos and bounded in $R_2$.

Expanding $g(\bar{y}, u^{**})$ about the point $(\bar{Y}, 1)$ in a second order Taylor series, we have

$$g(\bar{y}, u^{**}) = g[\bar{Y} + (\bar{y} - \bar{Y}), 1 + (u^{**} - 1)]$$
$$= g(\bar{Y}, 1) + (\bar{y} - \bar{Y})G_0 + (u^{**} - 1)G_1 + (u^{**} - 1)^2 G_2 + (\bar{y} - \bar{Y})(u^{**} - 1)G_3 + (\bar{y} - \bar{Y})G_4$$
$$= \bar{y} + (u^{**} - 1)G_1 + (u^{**} - 1)^2 G_2 + (\bar{y} - \bar{Y})(u^{**} - 1)G_3 + (\bar{y} - \bar{Y})G_4$$

$$\hat{\bar{Y}}_1 = \bar{Y} + \kappa_Y - \dfrac{n_1 \kappa_X}{\mu_x} G_1 + \dfrac{n_1^2 \kappa_X^2}{\mu_x^2} G_2 - \dfrac{n_1 \kappa_Y \kappa_X}{\mu_x} G_3 + \kappa_Y^2 G_4$$

$$\hat{\bar{Y}}_1 - \bar{Y} = \kappa_Y - \dfrac{n_1 \kappa_X}{\mu_x} G_1 + \dfrac{n_1^2 \kappa_X^2}{\mu_x^2} G_2 - \dfrac{n_1 \kappa_Y \kappa_X}{\mu_x} G_3 + \kappa_Y^2 G_4 \tag{26}$$

Taking expectations on both sides of (26) and using the definition of bias, we obtain

$$B(\hat{\bar{Y}}_1) = \dfrac{n_1^2 r_1}{\mu_x^2} G_2 - \dfrac{n_1 r_{01}}{\mu_x} G_3 + r_0 G_4 \tag{27}$$

By the definition of mean square error, we have

$$MSE(\hat{\bar{Y}}_1) = E[\hat{\bar{Y}}_1 - \bar{Y}]^2 = \left[\kappa_Y - \dfrac{n_1 \kappa_X}{\mu_x} G_1 + O(\kappa) - \bar{Y}\right]^2$$

$$= \kappa_Y^2 + \dfrac{n_1^2 \kappa_X^2}{\mu_x^2} G_1^2 - \dfrac{2 n_1 \kappa_Y \kappa_X}{\mu_x} G_1$$

$$\text{MSE}(\hat{\bar{Y}}_1) = r_0 + \frac{n_1^2 r_1}{\mu_x^2} G_1^2 - \frac{2n_1 r_{01}}{\mu_x} G_1 \tag{28}$$

On differentiating (28) with respect to $G_1$ and equating to zero we obtain

$$G_1(\text{opt}) = \frac{r_{01} \mu_x}{n_1 r_1} = G_1^{\Theta} \text{ (say)} \tag{29}$$

Using (28) and (29), we have the minimum MSE of $\hat{\bar{Y}}_1$ as

$$\min \text{MSE}(\hat{\bar{Y}}_1) = \left[ r_0 - \frac{r_{01}^2}{r_1} \right] \tag{30}$$

Using (25), we have the following particular members of $\hat{\bar{Y}}_1$ as

$$\hat{\bar{Y}}_1^1 = \bar{y}\varepsilon_1 + (1-\varepsilon_1)\frac{\bar{x}^{**}}{\mu_x} = \bar{y}\left[\varepsilon_1 + (1-\varepsilon_1)u^{**}\right] \tag{31}$$

$$\hat{\bar{Y}}_1^2 = \bar{y}\left[2 - \left(\frac{\mu_x}{\bar{x}^{**}}\right)^{\varepsilon_2}\right] = \bar{y}\left[2 - (u^{**})^{-\varepsilon_2}\right] \tag{32}$$

$$\hat{\bar{Y}}_1^3 = \bar{y}\left[\frac{\mu_x + \varepsilon_3(\bar{x}^{**} - \mu_x)}{\mu_x}\right] = \bar{y}\left[1 + \varepsilon_3(u^{**} - 1)\right] \tag{33}$$

$$\hat{\bar{Y}}_1^4 = \bar{y}\left[\frac{\mu_x + \varepsilon_3(\bar{x}^{**} - \mu_x)}{\bar{x}^{**}}\right] = \bar{y}\left[(u^{**})^{-1} + \varepsilon_3\left\{1 - (u^{**})^{-1}\right\}\right] \tag{34}$$

### 3.2.2 A Modified difference class of estimator

Motivated by Shabbir and Yaab (2003), we suggested a modified class of estimator given as

$$\hat{\bar{Y}}_2 = (1-J)\bar{y} + J\frac{t_b}{\mu_x} \tag{35}$$

where $t_b = \bar{y}\bar{x}^{**}\left(\frac{1+\gamma C_{yx}}{1+\gamma C_x^2}\right) \approx \bar{y}\bar{x}^{**}\lambda$ (say) and J is constant to be optimise.

Expressing $\hat{\bar{Y}}_2$ in terms of $\kappa_i$'s, we have

$$\hat{\bar{Y}}_2 = \left[J\lambda\left\{\bar{Y} + \kappa_Y - n_1 R\kappa_X - \frac{n_1 \kappa_Y \kappa_X}{\mu_x}\right\} + (1-J)\left\{\bar{Y} + \kappa_Y\right\}\right] \tag{36}$$

where $R = \frac{\bar{Y}}{\mu_x}$

Subtracting $\bar{Y}$ from both sides of equation (36) and then taking expectations, we have

$$B(\hat{\bar{Y}}_2) = \left[J\lambda\left\{\bar{Y} - \frac{n_1 r_{01}}{\mu_x}\right\} + (1-J)\bar{Y} - \bar{Y}\right] \tag{37}$$

By the definition of mean square error, we have

$$\text{MSE}(\hat{\bar{Y}}_2) = E\left[\hat{\bar{Y}}_2 - \bar{Y}\right]^2 = \left[J\lambda\left\{\bar{Y} + \kappa_Y - n_1R\kappa_X - \frac{n_1\kappa_Y\kappa_X}{\mu_x}\right\} + (1-J)\left\{\bar{Y} + \kappa_Y\right\}\right]^2 \quad (38)$$

$$= E\left[\bar{Y}^2 + J^2\lambda^2\left\{\bar{Y}^2 + \kappa_Y^2 + n_1^2R^2\kappa_X^2 - 4n_1R\kappa_Y\kappa_X\right\} + (1-J)^2\left\{\bar{Y}^2 + \kappa_Y^2\right\} - 2J\lambda\left\{\bar{Y}^2 - n_1R\kappa_Y\kappa_X\right\} \right.$$
$$\left. - 2(1-J) + 2J(1-J)\lambda\left\{\bar{Y}^2 + \kappa_Y^2 - 2n_1R\kappa_Y\kappa_X\right\}\right] \quad (39)$$

**Remark:** In the above equation (39), we have considered the terms up to the first order of approximation and neglecting terms whose expected value is assumed to be zero. Thus, we have

$$\text{MSE}(\hat{\bar{Y}}_2) = \left[\bar{Y}^2 + J^2C_1 + (1-J)^2C_2 - 2JC_3 - 2(1-J)C_4 + 2J(1-J)C_5\right] \quad (40)$$

The $\text{MSE}(\hat{\bar{Y}}_2)$ at (40) is minimised for

$$J(\text{opt}) = \frac{C_2 + C_3 - C_4 - C_5}{C_1 + C_2 - 2C_5} = J^\Theta \text{ (say)}$$

Thus the resulting minimum MSE of $\hat{\bar{Y}}_2$ is given by

$$\min \text{MSE}(\hat{\bar{Y}}_2) = \left[\left(\bar{Y}^2 + C_2 - 2C_4\right) - \frac{(C_2 + C_3 - C_4 - C_5)^2}{(C_1 + C_2 - 2C_5)}\right] = \left[\bar{Y}^2 + \phi_2\right] \quad (41)$$

where $\phi_2 = (C_2 - 2C_4) - \dfrac{(C_2 + C_3 - C_4 - C_5)^2}{(C_1 + C_2 - 2C_5)}$ and

$C_1 = \lambda^2\left[\bar{Y}^2 + r_0 + n_1^2R^2r_1 - 4n_1Rr_{01}\right]$, $C_2 = \lambda^2\left[\bar{Y}^2 + r_0\right]$, $C_3 = \lambda\left[\bar{Y}^2 - n_1Rr_{01}\right]$, $C_4 = \bar{Y}^2$
and $C_5 = \lambda\left[\bar{Y}^2 + r_0 - 2n_1Rr_{01}\right]$.

### 3.2.3 Adapted difference cum-dual-to ratio type estimator

Motivated by Singh and Solanki (2012), we propose a difference cum dual-to-ratio type estimator as

$$\hat{\bar{Y}}_P = d_1\bar{y}_\beta^* + d_2\bar{y}\left[\frac{c_1\bar{x}^{**} + c_2}{c_1\mu_x + c_2}\right]^{c_3} \quad (42)$$

where $\bar{y}_\beta^* = \bar{y} + \beta(\mu_x - \bar{x}^{**})$ is usual regression estimator, $(d_1, d_2)$ are suitably chosen scalars, $(c_1, c_2)$ are either constants or function of some known population parameter such as population mean $\mu_x$, population mean square $S_x^2$, correlation coefficient of variation $C_x$ and correlation coefficient between y and x ($\rho$). Also $c_3$ takes values (0, 1, -1) in order to make different ratio and product type estimators. Further, some particular members of $\hat{\bar{Y}}_P$ are listed in Table A.1 in appendix.

**Note** Here in equation (42) $\beta = \dfrac{S_{yx}}{S_x^2}$ is regression coefficient, which is assumed to be known.

Expressing (42) in terms of $\kappa_i$'s, we have

$$\hat{\bar{Y}}_P = d_1\left[\bar{Y} + \kappa_Y + \beta n_1\kappa_X\right] + d_2\bar{y}\left[1 - \tau_i n_1\kappa_X\right]^{c_3} \quad (43)$$

where $\tau_i = (c_1/c_1\mu_x + c_2)$ contribute the following possible values under which adapted estimator performs better, are given as

$$\tau_1 = [\rho/\rho\mu_x - C_x], \ \tau_2 = [1/\mu_x - C_x^2], \quad \tau_3 = [\rho/\rho\mu_x + C_x], \quad \tau_4 = [\rho/\rho\mu_x - C_x],$$
$$\tau_5 = [C_x/\mu_x(C_x - 1)] \ \tau_6 = [C_x/\mu_x(C_x + 1)], \ \tau_7 = [1/\mu_x + C_x] \text{ and } \tau_8 = [1/\mu_x - C_x]$$

We assume that $|\tau_i n_1 \kappa_X| < 1$, so that the term $[1 - \tau_i n_1 \kappa_X]^{c_3}$ is expandable. Thus by expanding the right hand side of (43) and neglecting the terms of $\kappa_i$'s having power greater than two, we have

$$\hat{\bar{Y}}_{Pi} - \bar{Y} = \left| d_1 [\bar{Y} + \kappa_Y + \beta n_1 \kappa_X] + d_2 \left[ \bar{Y} + \kappa_Y - c_3 \tau_i n_1 \kappa_X \bar{Y} - c_3 \tau_i n_1 \kappa_Y \kappa_X + \frac{c_3(c_3 - 1)}{2} \tau_i^2 n_1^2 \kappa_X^2 \bar{Y} \right] \right.$$
$$\left. - \bar{Y} \right| \quad (44)$$

Taking expectations on both sides of equation (44), we have

$$B(\hat{\bar{Y}}_{Pi}) = d_1 \bar{Y} + d_2 \left[ \bar{Y} - c_3 \tau_i n_1 r_{01} + \frac{c_3(c_3 - 1)}{2} \tau_i^2 n_1^2 \kappa_X^2 \bar{Y} \right] - \bar{Y} \quad (45)$$

By the definition of mean square error, we have
$$\text{MSE}(\hat{\bar{Y}}_P) = E[\hat{\bar{Y}}_P - \bar{Y}]^2$$

$$E(\hat{\bar{Y}}_{Pi} - \bar{Y})^2 = E[\bar{Y}^2 + d_1^2 \{\bar{Y}^2 + \kappa_Y^2 + \beta^2 n_1^2 \kappa_X^2 + 2\beta n_1 \kappa_Y \kappa_X\} + d_2^2 \{\bar{Y}^2 + \kappa_Y^2 + c_3^2 \tau_i^2 n_1^2 \kappa_X^2 \bar{Y}^2$$

$$- 4c_3 \tau_i n_1 \kappa_Y \kappa_X \bar{Y} + c_3(c_3 - 1)\tau_i^2 n_1^2 \kappa_X^2 \bar{Y}^2\} - 2d_1 \bar{Y}^2 - 2d_2 \left\{ \bar{Y} - c_3 \tau_i n_1 \kappa_Y \kappa_X + \frac{c_3(c_3 - 1)}{2} \tau_i^2 n_1^2 \kappa_X^2 \bar{Y} \right\}$$

$$+ 2d_1 d_2 \left\{ \bar{Y}^2 + \kappa_Y^2 - 2c_3 \tau_i n_1 \kappa_Y \kappa_X + \frac{c_3(c_3 - 1)}{2} \tau_i^2 n_1^2 \kappa_X^2 \bar{Y}^2 + \beta n_1 \kappa_Y \kappa_X - c_3 \beta \tau_i n_1^2 \kappa_X^2 \bar{Y} \right\}$$

$$\text{MSE}(\hat{\bar{Y}}_P) = [\bar{Y}^2 + d_1^2 D_1 + d_2^2 D_2 - 2d_1 D_3 - 2d_2 D_4 + 2d_1 d_2 B_5] \quad (46)$$

The $\text{MSE}(\hat{\bar{Y}}_P)$ at (46) is minimised for

$$d_1(\text{opt}) = \left( \frac{D_2 D_3 - D_4 D_5}{D_1 D_2 - D_5^2} \right) = d_1^\Theta \text{ (say)}$$

$$d_2(\text{opt}) = \left( \frac{D_1 D_4 - D_3 D_5}{D_1 D_2 - D_5^2} \right) = d_2^\Theta \text{ (say)}$$

Thus the resulting minimum MSE of $(\hat{\bar{Y}}_P)$ is given by

$$\text{MSE}(\hat{\bar{Y}}_P) = \left[ \bar{Y}^2 - \frac{(D_1 D_2^2 D_3^2 - D_1 D_4^2 D_5^2 + D_1^2 D_2 D_4^2 - D_2 D_3^2 D_5^2 + 2D_3 D_4 D_5^3 - 2D_1 D_2 D_3 D_4 D_5)}{(D_1 D_2 - D_5^2)^2} \right]$$

$$\min \text{MSE}(\hat{\bar{Y}}_P) = \left[\bar{Y}^2 - \phi_P\right] \tag{47}$$

where $\phi_P = \dfrac{\left(D_1 D_2^2 D_3^2 - D_1 D_4^2 D_5^2 + D_1^2 D_2 D_4^2 - D_2 D_3^2 D_5^2 + 2 D_3 D_4 D_5^3 - 2 D_1 D_2 D_3 D_4 D_5\right)}{\left(D_1 D_2 - D_5^2\right)^2}$.

We would like to mention here that the proposed class of estimator $\hat{\bar{Y}}_P$ is reduced to some known estimators of $\bar{Y}$ by putting different values of $\left([d_1, d_2, c_1, c_2, c_3]\right)$ ie.

$[d_1, d_2, c_1, c_2, c_3] = [0,1,1,0,1]; \hat{\bar{Y}}_P \to e_1^m \to$ Srivenkataramana (1980) estimator,

$[d_1, d_2, c_1, c_2, c_3] = [0,1,1,0,-1]; \hat{\bar{Y}}_P \to$ Dual-to-product type estimator,

$[d_1, d_2, c_1, c_2, c_3] = [1,0,-,-,-]; \hat{\bar{Y}}_P \to$ Usual regression estimator.

## 4 Efficiency Comparisons

From (12), (17), (17), (23), (30), (41) and (47) we have

$\text{MSE}(\hat{\bar{Y}}_1) < V(\bar{y})$, if

$$\left[r_0 - \frac{r_{01}^2}{r_1}\right] - \gamma \bar{Y}^2 \left(C_Y^2 + \frac{S_{d_Y}^2}{\bar{Y}^2}\right) \leq 0 \tag{48}$$

$\text{MSE}(\hat{\bar{Y}}_2) < V(\bar{y})$, if

$$\gamma \bar{Y}^2 \left(C_Y^2 + \frac{S_{d_Y}^2}{\bar{Y}^2}\right) - \left(\bar{Y}^2 + \phi_2\right) \geq 0 \tag{49}$$

$\text{MSE}(\hat{\bar{Y}}_1) < \text{MSE}(e_1^m)$, if

$$\gamma \bar{Y}^2 \left[C_Y^2 + n_1^2 C_X^2 - 2 n_1 \rho C_Y C_X\right] + \gamma \left[S_{d_Y}^2 + n_1^2 R^2 S_{d_X}^2\right] - \left[r_0 - \frac{r_{01}^2}{r_1}\right] \geq 0 \tag{50}$$

$\text{MSE}(\hat{\bar{Y}}_2) < \text{MSE}(e_1^m)$, if

$$\gamma \bar{Y}^2 \left[C_Y^2 + n_1^2 C_X^2 - 2 n_1 \rho C_Y C_X\right] + \gamma \left[S_{d_Y}^2 + n_1^2 R^2 S_{d_X}^2\right] - \left(\bar{Y}^2 + \phi_2\right) \geq 0 \tag{51}$$

$\text{MSE}(\hat{\bar{Y}}_P) < V(\bar{y})$, if

$$\left[\bar{Y}^2 - \phi_P\right] - \gamma \bar{Y}^2 \left(C_Y^2 + \frac{S_{d_Y}^2}{\bar{Y}^2}\right) \leq 0 \tag{52}$$

$\text{MSE}(\hat{\bar{Y}}_P) < \text{MSE}(e_1^m)$, if

$$\gamma \bar{Y}^2 \left[C_Y^2 + n_1^2 C_X^2 - 2 n_1 \rho C_Y C_X\right] + \gamma \left[S_{d_Y}^2 + n_1^2 R^2 S_{d_X}^2\right] - \left(\bar{Y}^2 + \phi_P\right) \geq 0 \tag{53}$$

$$\mathrm{MSE}(\hat{\bar{Y}}_P) < \mathrm{MSE}(e_2^m), \text{ if}$$

$$\left[\bar{Y}^2 + \alpha'^2 B_1 + (1-\alpha')^2 B_2 - 2\alpha' B_3 - 2(1-\alpha') B_4 + 2\alpha(1-\alpha') B_5\right] - \left(\bar{Y}^2 + \phi_P\right) \geq 0 \quad (54)$$

If the above condition (48-54) holds, adapted class $\left[\hat{\bar{Y}}_1, \hat{\bar{Y}}_2, \hat{\bar{Y}}_P\right]$ performs much better than existing one.

## 5. Empirical Study

To evaluate the performance of adapted estimators $\left(\hat{\bar{Y}}_1, \hat{\bar{Y}}_2, \hat{\bar{Y}}_{Pi}\right)$ over other competitors, we have considered two population data sets for sample size n=500. The description of these data sets is as follows.

**Population 1**

$X = N(5,10), Y = X + N(0,1), y = Y + N(1,3), x = X + N(1,3), N=5000, \bar{Y} = 4.927167$

$\bar{X} = 4.924306, S_Y^2 = 102.0075, S_X^2 = 101.4117, S_{d_y}^2 = 8.862114, S_{d_x}^2 = 24.19283, \rho = 0.995059$

**Population 2**

$X = N(5,10), Y = X + N(0,1), y = Y + N(1,5), x = X + N(1,5), N=5000, \bar{Y} = 4.996681$

$\bar{X} = 5.013507, S_Y^2 = 97.12064, S_X^2 = 95.95803, S_{d_y}^2 = 23.96055, S_{d_x}^2 = 24.19283, \rho = 0.994822$

We have computed the percent relative efficiencies (PREs) of different estimators T, with respect to usual unbiased estimator $\bar{y}$ as

$$\mathrm{PREs}\,(T, \bar{y}) = \frac{\mathrm{Var}(\bar{y})}{\mathrm{MSE}_{\min}(T)} * 100$$

And the result are displayed in Table 1

**Table 1** Shows PREs and MSE's of adapted and existing estimators considered in section 3.1.

| Estimators | Population I | Population II |
|---|---|---|
| | PRE/MSE | PRE/MSE |
| $\bar{y}$ | 100/0.19956 | 100/ 0.217946 |
| $e_1$ | 123.56/0.16151 | 119.55/0.182305 |
| $e_2$ | 612.48/0.03258 | 273.214/0.079771 |
| $\hat{\bar{Y}}_1$ | 612.48/0.03258 | 273.214/0.079771 |
| $\hat{\bar{Y}}_2$ | 611.66/0.03263 | 273.2932/0.079748 |
| $\hat{\bar{Y}}_P^1$ | 618.29/ 0.032276 | 273.2585/0.079758 |
| $\hat{\bar{Y}}_P^2$ | 940.53/0.021218 | 315.404/0.069101 |
| $\hat{\bar{Y}}_P^3$ | 959.49/0.020799 | 302.231/0.072112 |
| $\hat{\bar{Y}}_P^4$ | 834.3038/0.02392 | 288.736/0.075483 |
| $\hat{\bar{Y}}_P^5$ | 822.301/0.024269 | 298.442/0.073028 |
| $\hat{\bar{Y}}_P^6$ | 945.54/0.021106 | **315.8539/0.069** |
| $\hat{\bar{Y}}_P^7$ | **964.96/0.020681** | 302.6126/0.072021 |

From Table 1 we conclude that adapted classes $(\hat{\bar{Y}}_1, \hat{\bar{Y}}_2)$ are better than usual unbiased estimator $\bar{y}$ and Srivenkataramana estimator $e_2$. Further, the proposed class of estimators $\hat{\bar{Y}}_P$ which utilizes the information on several population parameters of auxiliary variable x has an improvement over regression method of estimation and other existing estimators of population mean $\bar{Y}$ which utilizes the information only on population mean of auxiliary variable x. Among all, $\hat{\bar{Y}}_P^7$ is the best one for application point of view.

## 6. Conclusion

In this article we have suggested three different classes of estimators for estimating population mean $\bar{Y}$ in the presence of measurement error. The asymptotic bias and mean square error formulae of proposed classes have been obtained. The asymptotic optimum estimators in the proposed classes have been identified with its properties. It has been identified theoretically and numerically in section 4 and section 5 the proposed class $\hat{\bar{Y}}_P$ is better then all the estimators considered in section 3.1. Thus the proposed class $\left[\hat{\bar{Y}}_1, \hat{\bar{Y}}_2, \hat{\bar{Y}}_P\right]$ of estimators has been recommended for its use in practice.

## 7. References

[1] Allen, J., Singh, H.P. and Smarandache, F. (2003): A family of estimators of population mean using multi-auxiliary information in presence of measurement errors. Int. J. Soc. Econ., 30(7): 837-849.
[2] Bahl, S. and Tuteja, R.K. (1991): Ratio and product type exponential estimator. Infrm. and Optim. Sci., XII, I: 159-163.
[3] Cochran, W.G. (1968): Errors of measurement in statistics. Technometrics, 10(4): 637-666.

[4] Cochran, W. G. (2005): Sampling techniques. 2nd Ed. New Delhi, India: Wiley Eastern Private Limited.

[5] Kumar, M., Singh, R., Singh, A.K. and Smarandache, F. (2011): Some ratio type estimators under measurement errors. World Applied Sciences Journal, 14(2): 272-276.

[6] Manisha and Singh, R.K. (2001): An estimation of population mean in the presence of measurement errors. J. Ind. Soc. Agri. Statist., 54(1): 13-18.

[7] Manisha and Singh, R.K. (2002): Role of regression estimator involving measurement errors. Brazilian J. Probability Statistics, 16: 39-46.

[8] Shalabh (1997): Ratio method of estimation in the presence of measurement errors. J. Ind. Soc. Agri. Statist., 50(2): 150-155.

[9] Sharma B. and Tailor R. (2010): A new ratio-cum-dual to ratio estimator of finite population mean in simple random sampling, Global Journal of Science Frontier Research, 10(1), 27-31.

[10] Shabbir, J., Yaab, M. Z. (2003): Improvement over transformed auxiliary variable in estimating the finite population mean. *Biometrical J.* 45:723–729.

[11] Shukla, D., Pathak, S. and Thakur, N. S.(2012): An estimator for mean estimation in presence measurement error. *Research & Reviews: A Journal of Statistics*, *1*(*1*), 1-8.

[12] Singh, H.P. and Solanki, R.S. (2012): Improved estimation of population mean in simple random sampling using information on auxiliary attribute. Appl. Math. and Comp., 218: 7798-7812.

[13] Singh, H. P., & Karpe, N. (2008): Ratio- Product estimator for population mean in presence of measurement errors. Journal of Applied Statistical Sciences, 16(4), 49-64.

[14] Singh, H.P. and Karpe, N. (2009): On estimation of two populations means using supplementary information in presence of measurement errors. Statistica, Department of Statistics, University of Bologna, 69(1): 27-47.

[15] Singh, V. K., **Singh, R.** and Florentin, S. (2014): Difference-type estimators for estimation of mean in the presence of measurement error. *The efficient use of supplementary information in finite population sampling. Education Publishing*, USA**.** ISBN : 978-1-59973-275-6.

[16] Srivastava, A.K. and Shalabh (2001): Effect of measurement errors on the regression method of estimation in survey sampling. J. Statist. Res., 35(2): 35-44.

[17] Srivastava, S. K. (1971): A generalized estimator for the mean of a finite population using multi auxiliary information, Jour. Amer. Stat. Assoc., 66, 404-407.

[18] Srivenkataramana, T. (1980): A dual to ratio estimator in sample surveys. Biometrika, 67(1): 199-204.

[19] Sukhatme, P. V., Sukhatme, B. V., Sukhatme, S., & Ashok, C. (1984): Sampling theory of surveys with applications. New Delhi: Iowa State University Press.

**Appendix**

In table A.1 listed below have some members of proposed class of estimators $\hat{\bar{Y}}_P$ given as

Some particular members of proposed class $\hat{\bar{Y}}_P$

| Estimator | Different parameters | | |
|---|---|---|---|
| | $c_1$ | $c_2$ | $c_3$ |
| $\hat{\bar{Y}}_P^1 = d_1 \bar{y}_\beta^* + d_2 \bar{y} \left[ \dfrac{C_x - \rho \bar{x}^{**}}{C_x - \rho \mu_x} \right]$ | $-\rho$ | $C_x$ | 1 |
| $\hat{\bar{Y}}_P^2 = d_1 \bar{y}_\beta^* + d_2 \bar{y} \left[ \dfrac{C_x + \rho \bar{x}^{**}}{C_x + \rho \mu_x} \right]^{-1}$ | $\rho$ | $C_x$ | -1 |
| $\hat{\bar{Y}}_P^3 = d_1 \bar{y}_\beta^* + d_2 \bar{y} \left[ \dfrac{C_x - \rho \bar{x}^{**}}{C_x - \rho \mu_x} \right]^{-1}$ | $-\rho$ | $C_x$ | -1 |
| $\hat{\bar{Y}}_P^4 = d_1 \bar{y}_\beta^* + d_2 \bar{y} \left[ \dfrac{\bar{X} - C_x \bar{x}^{**}}{\bar{X} - C_x \mu_x} \right]^{-1}$ | $-C_x$ | $\bar{X}$ | -1 |
| $\hat{\bar{Y}}_P^5 = d_1 \bar{y}_\beta^* + d_2 \bar{y} \left[ \dfrac{\bar{X} + C_x \bar{x}^{**}}{\bar{X} + C_x \mu_x} \right]^{-1}$ | $C_x$ | $\bar{X}$ | -1 |
| $\hat{\bar{Y}}_P^6 = d_1 \bar{y}_\beta^* + d_2 \bar{y} \left[ \dfrac{C_x + \bar{x}^{**}}{C_x + \mu_x} \right]^{-1}$ | 1 | $C_x$ | -1 |
| $\hat{\bar{Y}}_P^7 = d_1 \bar{y}_\beta^* + d_2 \bar{y} \left[ \dfrac{\bar{x}^{**} - C_x}{\mu_x - C_x} \right]^{-1}$ | 1 | $-C_x$ | -1 |